\documentclass[onefignum,onetabnum]{siamart171218}

\usepackage{amsmath}
\usepackage{amssymb}
\usepackage{blkarray}
\usepackage{bm}
\usepackage{color}
\usepackage{tabularx}
\usepackage{caption}
\usepackage{multirow}
\usepackage{graphicx}

\usepackage{algpseudocode}

\usepackage{subcaption}
\usepackage{mathtools}
\usepackage{stmaryrd}
\usepackage{url}

\usepackage{subcaption}

\newsiamremark{remark}{Remark}
\newsiamremark{hypothesis}{Hypothesis}
\crefname{hypothesis}{Hypothesis}{Hypotheses}
\newsiamthm{claim}{Claim}

\newcommand{\TheTitle}{A Scalable Multigrid Reduction Framework for Multiphase Poromechanics of Heterogeneous Media}
\newcommand{\ShortTitle}{Multigrid Reduction Framework for Multiphase Poromechanics}
\newcommand{\TheAuthors}{Q. M. Bui, D. Osei-Kuffuor, N. Castelletto, J. A. White}

\headers{\ShortTitle}{\TheAuthors}

\title{\TheTitle\footnotemark[1]}

\author{Quan M. Bui\footnotemark[2]
\and Daniel Osei-Kuffuor\footnotemark[3]
\and Nicola Castelletto\footnotemark[4]
\and Joshua A. White\footnotemark[4]}

\begin{document}
\maketitle

\renewcommand{\thefootnote}{\fnsymbol{footnote}}
\footnotetext[1]{Submitted to the editors \today.
\funding{This work was funded by LLNL through Laboratory Directed Research and Development (LDRD) Project 18-ERD-027.  Portions of this work were performed under the auspices of the U.S. Department of Energy by Lawrence Livermore National Laboratory under Contract DE-AC52-07-NA27344.}}
\footnotetext[2]{Corresponding Author. Center for Applied Scientific Computing, Lawrence Livermore National Laboratory, United States (\email{bui9@llnl.gov})}
\footnotetext[3]{Center for Applied Scientific Computing, Lawrence Livermore National Laboratory, United States (\email{oseikuffuor1@llnl.gov}).}
\footnotetext[4]{Atmospheric, Earth and Energy Division, Lawrence Livermore National Laboratory, United States (\email{castelletto1@llnl.gov}, \email{jawhite@llnl.gov}).}

\begin{abstract}
Simulation of multiphase poromechanics involves solving a multi-physics problem in which multiphase flow and transport are tightly coupled with the porous medium deformation.
To capture this dynamic interplay, fully implicit methods, also known as monolithic approaches, 
are usually preferred.
The main bottleneck of a monolithic approach is that it requires solution of large linear systems that result from the discretization and linearization of the governing balance equations.
Because such systems are non-symmetric, indefinite, and highly ill-conditioned, preconditioning is 
critical for fast convergence.
Recently, most efforts in designing efficient preconditioners for multiphase poromechanics have been dominated by physics-based strategies.
Current state-of-the-art ``black-box'' solvers such as algebraic multigrid (AMG) are ineffective because they cannot effectively capture the strong coupling between the mechanics and the flow sub-problems, as well as the coupling inherent in the multiphase flow and transport process. 
In this work, we develop an algebraic framework based on multigrid reduction (MGR) that is suited for tightly coupled systems of PDEs.
Using this framework, the decoupling between the equations is done algebraically through defining appropriate interpolation and restriction operators.
One can then employ existing solvers for each of the decoupled blocks or design a new solver based on knowledge of the physics.
We demonstrate the applicability of our framework when used as a ``black-box'' solver for multiphase poromechanics.
We show that the framework is flexible to accommodate a wide range of scenarios, as well as efficient and scalable for large problems.
\end{abstract}

\section{Introduction} \label{sec:intro}
Modeling subsurface systems requires an understanding of many different physical processes, including multiphase fluid flow and transport, and geomechanical deformations.
These processes are often tightly coupled in a ``two-way'' fashion: for example, the flow process can have a large influence on the mechanical process and in turn can be affected by the feedback of the produced mechanical response.
To simulate these processes, one needs to solve a set of coupled, nonlinear, time-dependent partial differential equations (PDEs) that govern the conservation of mass of the fluid phases and linear momentum of the solid-fluid mixture.
For this system, fully-implicit time discretization is the widely preferred approach because it is unconditionally stable and allows for large time steps.
However, using an implicit approach, one must solve a large, sparse, and ill-conditioned linear system at each nonlinear iteration.
Robust and scalable solvers are therefore needed for large scale simulations on high performance computing platforms.
This paper presents our efforts to design an efficient preconditioning strategy based on an algebraic framework that is flexible and capable of addressing the inherent ill-conditioning of a complicated multi-physics system.

In recent years, much of the work in developing preconditioning strategies for coupled problems has focused on so-called physics-based strategies.
The key is to use knowledge of the specific physical processes involved to break the tightly coupled systems into smaller sub-problems whose properties are well-studied.
For example, these sub-problems could take the form of an elliptic, hyperbolic, or parabolic PDE, to which appropriate techniques can be applied.
For fully implicit simulation of complex multiphase flow and transport without mechanics, one of the most popular methods is the \textit{Constrained Residual Pressure} (CPR) multistage preconditioning technique \cite{Wallis83,Wallis85}.
For single-phase flow poromechanics, many block preconditioners have been developed \cite{Adler18,Bergamaschi07,Bergamaschi12,Haga12,Lee17,White11,White16}.
In the context of multiphase poromechanics, one recent strategy \cite{White18} uses the \textit{fixed-stress} partitioning \cite{KimSPE11,Mikeli12,Settari98,White16} of the mechanics and the flow parts combined with a CPR approach \cite{Wallis83} for the multiphase flow system.
In general, these physics-based preconditioners are among the most effective techniques available.
However, designing a good strategy is both time-consuming and challenging as it requires extensive knowledge of the particular continuous model of interest.
One alternative is to use a ``black-box'' approach, such as \textit{algebraic multigrid} (AMG) \cite{Stueben01,Trottenberg00}.
AMG techniques are among the most efficient and scalable methods for solving sparse linear systems.
Unlike \textit{geometric multigrid}, these methods do not need an explicit hierarchy of computational grids.
However, they are originally designed for scalar elliptic PDEs, and their applicability is much more limited for strongly coupled systems of PDEs.

Recently, \textit{multigrid reduction} (MGR) \cite{Ries79,Ries83}, a variant of AMG, has been applied successfully to coupled systems of multiphase flow and transport with phase transitions \cite{Bui18,Wang17}.
Drawing on this success, in this work we further develop MGR into a general multi-level framework for solving discrete systems coming from discretization of tightly coupled PDEs.
We also introduce a new dropping strategy for computing the reduction onto the coarse-grid within the MGR V-cycle that effectively captures the coupling between mechanics and flow.
The goal of this strategy is two-pronged: (1) to keep the coarse grid sparse during multi-level reduction, and (2) to make the coarse grid amenable to classical AMG.
We show that with this new feature, MGR is effective as a general-purpose algebraic solver for multiphase 
poromechanics, and it also scales well with problem size.
The rest of the paper is organized as follows.
Section 2 and 3 introduce the governing equations of multiphase poromechanics and the discretization scheme.
Section 4 describes the nonlinear solution algorithm.
In section 5, we describe the MGR framework and how it is applied to solve the linear systems coming from the linearization.
We show numerical results in section 6 to demonstrate the robustness and scability of the proposed preconditioner.
We then end with some concluding remarks and directions for future work.


\section{Problem Statement} \label{sec:statement}
We focus on a displacement-saturation-pressure formulation of immiscible two-phase flow through a deforming poroelastic medium \cite{Cou04}.
We limit the discussion to quasi-static small-strain kinematics.
Let the subscript $w$ and $nw$ denote the wetting and non-wetting fluid phase, respectively. 
Since the medium's pore space is always fluid-filled, the fluid phase saturations must always sum to unity, i.e. $(s_w + s_{nw}) = 1$.
Here, the wetting fluid phase saturation, denoted from now on by lower case $s$ without subscript, is used as a primary unknown.
Capillary pressure, which is the difference between the phase pressure of the non-wetting phase and the wetting phase, is not considered---a frequent assumption in many practical engineering applications.
Hence, we have $p_w = p_{nw} = p$.

For a given closed domain $\overline{\Omega} = \Omega \cup \Gamma \in \mathbb{R}^3$, with $\Omega$ an open set and $\Gamma$ its boundary, and time interval $\mathbb{I} = (0, t_{\max}]$, the strong form of the multiphase poromechanical initial/boundary value problem (IBVP) consists of finding the displacement vector field $\boldsymbol{u} : \overline{\Omega} \times \mathbb{I} \rightarrow \mathbb{R}^3$, the wetting fluid phase saturation $s : \overline{\Omega} \times \mathbb{I} \rightarrow \mathbb{R}^3$, and the pore pressure $p: \overline{\Omega} \times \mathbb{I} \rightarrow \mathbb{R}$ such that \cite{Cou04}:
\begin{subequations}\label{eq:IBVP_global}
\begin{align}
  &-\nabla \cdot \boldsymbol{\sigma} = \rho \boldsymbol{g} & &\mbox{ on } \Omega \times \mathbb{I} & & \mbox{(linear momentum balance),} \label{eq:momentumBalanceS}\\
  %
  &\dot{m}_{w} + \nabla \cdot \boldsymbol{w}_{w} = q_{w} & &\mbox{ on } \Omega \times \mathbb{I} & &\mbox{(wetting fluid phase mass balance),} \label{eq:massBalanceW_S}	 \\
  %
  &\dot{m}_{nw} + \nabla \cdot \boldsymbol{w}_{nw} = q_{nw} & &\mbox{ on } \Omega \times \mathbb{I} & &\mbox{(non-wetting fluid phase mass balance).} \label{eq:massBalanceO_S}	
\end{align}
\end{subequations}
where 
\begin{itemize}
  \item $\boldsymbol{\sigma} = \left( \textbf{\sffamily C} : \nabla^s \boldsymbol{u}  - b p  \boldsymbol{1} \right)$ is the total Cauchy stress tensor, with $\textbf{\sffamily C}$ the rank-4 elasticity tensor, $b$ the Biot coefficient, and $\boldsymbol{1}$ the rank-2 identity tensor;
  \item $\rho \boldsymbol{g}$ is a body force due to the self-weight of the multiphase mixture, with $\rho =  ((1-\phi) \rho_s + \phi \rho_w s + \phi \rho_{nw} (1-s))$ the density of the mixture, $\phi$ the porosity, $\rho_s$, $\rho_w$, $\rho_{nw}$ the density of the solid, the wetting, and the non-wetting fluid phase, respectively, and $\boldsymbol{g}$ the gravity vector;
  \item $m_w = (\phi \rho_w s)$ and $m_{nw} = (\phi \rho_{nw} (1-s))$ denote wetting and non-wetting fluid phase mass per unit volume;
  \item $\boldsymbol{w}_w = - (\rho_w \lambda_w \boldsymbol{\kappa} \cdot \nabla \Phi_{w})$ and $\boldsymbol{w}_{nw} = - (\rho_{nw} \lambda_{nw} \boldsymbol{\kappa} \cdot \nabla \Phi_{nw})$ are wetting and non-wetting fluid phase mass fluxes \cite{Aziz79}, with $\lambda_{\ell} = k_{r{\ell}}/\mu_{\ell}$ the phase mobility, $\mu_{\ell}$ the phase viscosity, $k_{r\ell}$ the phase  relative permeability factor,  $\boldsymbol{\kappa}$ the absolute permeability tensor, $\Phi_{\ell} = (p - \rho_{\ell} \boldsymbol{g} \cdot \boldsymbol{x} )$ the phase potential, $\boldsymbol{x}$ the position vector in $\mathbb{R}^3$, $\ell = \{w, nw\}$;
  \item $q_w$ and $q_{nw}$ are mass source/sink per unit volume terms for the wetting and the non-wetting fluid phase, respectively;
  \item $\nabla$, $\nabla^s$, and $\nabla \cdot$ are the gradient, symmetric gradient, and divergence operator, respectively;
  \item the superposed dot, $\dot{(\bullet)}$, indicates the derivative of quantity $(\bullet)$ with respect to time .
\end{itemize}
For the application of boundary conditions, let us introduce two disjoint partitions of the domain boundary such that $\Gamma = \overline{\Gamma_u^D \cup \Gamma_u^N } = \overline{\Gamma_f^D \cup \Gamma_f^N }$.   
Without loss of generality, consider homogeneous displacement boundary conditions $\boldsymbol{u} = \boldsymbol{0}$ on $\Gamma_u^D \times \mathbb{I}$ and homogeneous flux conditions $\boldsymbol{w}_w \cdot \boldsymbol{n} = \boldsymbol{w}_{nw} \cdot \boldsymbol{n} = 0$ on $\Gamma_f^N \times \mathbb{I}$, along with a prescribed total traction conditions $\boldsymbol{\sigma} \cdot \boldsymbol{n} = \boldsymbol{t}^N$ on $\Gamma_u^N \times \mathbb{I}$ and pressure/saturation conditions $p = p^D$ and $s = s^D$ on $\Gamma_f^D \times \mathbb{I}$, where $\boldsymbol{n}$ denotes the outer normal vector for $\Gamma$. More complicated boundary conditions may be introduced as needed with modest changes to the discretization below.
The formulation is completed by appropriate: (i) initial conditions for $\boldsymbol{u}$, $s$, and $p$, and (ii) equations of state and constitutive equations to specify the following dependencies: $\phi = \phi(\boldsymbol{u}, p)$, $\rho_\ell = \rho_\ell(p)$, $\mu_\ell = \mu_\ell(p)$, and $k_{r\ell} = k_{r\ell} (s)$, with $\ell = \{ w, nw \}$.
For additional details on the adopted poromechanical model we refer the reader to \cite{White18}.

\section{Discretization}\label{sec:discretization}
Several space discretization methods for the multiphase po\-ro\-me\-cha\-ni\-cal IBVP have been proposed in the literature---see, e.g., \cite[and references therein]{LewSch98,HagOsnLan12b,White18}.
In this work, the discrete form of \eqref{eq:IBVP_global} is obtained by combining a finite element (FE) method for the mechanical subproblem with a finite volume (FV) approach for the multiphase flow and transport subproblem.
This choice is quite common when modeling nonlinear hydromechanical processes in subsurface formations characterized by highly heterogeneous hydrogeological properties, e.g. high-contrast permeability fields typically encountered in practical reservoir simulation \cite{Settari98,KimTchJua13,Pre14,GarKarTch16,Set_etal17}.

Let $\boldsymbol{H}_0^1(\Omega)$ denote the Sobolev space of vector functions satisfying displacement homogeneous Dirichlet conditions over $\Gamma_u^D$ and whose first derivatives belong to $L^2(\Omega)$, with $L^2(\Omega)$ the space of square integrable functions in $\Omega$; let $\boldsymbol{\mathcal{U}}^h \subset \mathbf{H}_0^1(\Omega)$, $\mathcal{S}^h \subset L^2(\Omega)$, $\mathcal{P}^h \subset L^2(\Omega)$ denote finite-dimensional subspaces associated with a conforming triangulation $\mathcal{T}^h$ of the domain into nonoverlapping hexahedral cells; and let $\widehat{w}_{\ell}^{\varepsilon}$ denote a conservative numerical flux approximating the $\ell$ fluid phase mass flux across face $\varepsilon$ in $\mathcal{E}^h$, namely the set of faces in $\mathcal{T}^h$, such that $\widehat{w}_{\ell}^{\varepsilon} \approx - \int_{\varepsilon} \boldsymbol{w}_\ell \cdot \boldsymbol{n}_{\varepsilon} \; \mathrm{d}A$, with $\boldsymbol{n}_{\varepsilon}$ a unit normal vector defining the global face orientation.

Precisely, our space discretization employs: (i) continuous piecewise trilinear finite elements for $\boldsymbol{\mathcal{U}}^h$, (ii) piecewise constant functions for $\mathcal{S}^h$ and $\mathcal{P}^h$, and (iii) a linear two-point flux approximation (TPFA) combined with a first-order upwinding strategy for the 
numerical flux $\widehat{w}_{\ell}^{\varepsilon}$ \cite{Aziz79}.
Using a fully-implicit time marching scheme, with the subscript $n$ indicating the discrete time level, the fully discrete mixed FE/FV variational statement of \eqref{eq:IBVP_global} is: find $\{\boldsymbol{u}^h_n, s^h_n, p^h_n\} \in \pmb{\mathcal{U}}^h \times \mathcal{S}^h \times \mathcal{P}^h$ such that for all $\{\boldsymbol{\eta}^h, \psi^h, \chi^h\} \in \pmb{\mathcal{U}}^h \times \mathcal{S}^h \times \mathcal{P}^h$
\begin{subequations}\label{eq:IBVP_G}
\begin{align}
   \mathcal{F}_u &= \left(\nabla^s \boldsymbol{\eta}^h,\textbf{\sffamily C} : \nabla^s \boldsymbol{u}_n^h \right) -
    \left(\nabla \cdot \boldsymbol{\eta}^h, bp_n^h \right)\ -
    \left(\boldsymbol{\eta}^h, \rho_n \boldsymbol{g} \right) -
    \int_{\Gamma_u^N} \boldsymbol{\eta}^h \cdot \boldsymbol{t}^N_n \; \mathrm{d}A = 0, \label{eq:IBVP_Gu} \\
   \mathcal{F}_s &= \left(\psi^h, \frac{m_{w,n} - m_{w,n-1}}{\Delta t_n} \right) - \sum_{\varepsilon \in \mathcal{E}^h \setminus \mathcal{E}_f^{h,N} } \llbracket \psi^h \rrbracket_{\varepsilon} \widehat{w}_{w,n}^{\varepsilon} - (\psi^h,q_{w,n}) = 0, \label{eq:IBVP_Gw} \\
   \mathcal{F}_p &= \left(\chi^h, \frac{m_{nw,n} - m_{nw,n-1}}{\Delta t_n} \right) - \sum_{\varepsilon \in \mathcal{E}^h \setminus \mathcal{E}_f^{h,N} } \llbracket \chi^h \rrbracket_{\varepsilon} \widehat{w}_{nw,n}^{\varepsilon} - (\chi^h,q_{nw,n}) = 0, \label{eq:IBVP_Gnw}
\end{align} 
\end{subequations}
where $n \in \{ 1, 2, \ldots \}$.
The compact notation $(\bullet,\bullet)$ denotes the $L^2$-inner product of scalar, vector, or rank-2 tensor functions in $L^2(\Omega)$, $[L^2(\Omega)]^3$, or $[L^2(\Omega)]^{3 \times 3}$, as appropriate.
In \eqref{eq:IBVP_Gw}-\eqref{eq:IBVP_Gnw}, $\Delta t_n = (t_n - t_{n-1})$ is the timestep size; $\mathcal{E}_f^{h,N}$ is the set of faces belonging to the boundary $\Gamma_f^{N}$; and $\llbracket \bullet \rrbracket_{\varepsilon}$ indicates the jump of a quantity $(\bullet)$ across $\varepsilon$.
For an internal face $\varepsilon$ shared by cells $K$ and $L$, $\boldsymbol{n}_{\varepsilon}$ pointing from $K$ to $L$, $\llbracket \bullet \rrbracket_{\varepsilon} = ( (\bullet)_{|L} - (\bullet)_{|K} )$, with $(\bullet)_{|K}$ and $(\bullet)_{|L}$ the restriction of $(\bullet)$ on $K$ and $L$, respectively.
For a boundary lying on $\Gamma_f^{D}$, $\boldsymbol{n}_{\varepsilon}$ coincides with the outer normal to the domain boundary and the jump expression simply reads $\llbracket \bullet \rrbracket_{\varepsilon} = - (\bullet)_{|K}$.

Finally, introducing in \eqref{eq:IBVP_G} the expressions $\boldsymbol{u}_n^h = \sum_i u_{i,n} \boldsymbol{\eta}_i^h$, $s_{n}^h = \sum_j s_{j,n} \psi_j^h$, and $p_{n}^h = \sum_k p_{k,n} \chi_k^h$, with $\{ \boldsymbol{\eta}_i^h \}$, $\{ \psi_j^h \}$, and $\{ \chi_k^h \}$ bases for $\boldsymbol{\mathcal{U}}^h$, $\mathcal{S}^h$, and $\mathcal{P}^h$, respectively,  a standard Galerkin approach yields a system of nonlinear discrete equations
\begin{align}
  F(\mathbf{x}_{n} ) &=
  \left(
  \begin{array}{c}
    F_u(\mathbf{x}_{n} ) \\
    F_s(\mathbf{x}_{n} ) \\
    F_p(\mathbf{x}_{n} ) 
  \end{array}
  \right)
  = 
  \mathbf{0}.
  \label{eq:IBVP_G_res}
\end{align}
Here, vector $\mathbf{x}_n$ contains the nodal displacement ($u_{i,n}$), cell-centered saturation $s_{i,n}$ and cell-centered pressure $p_{i,n}$ coefficients that are used to expand $\boldsymbol{u}^h_n$, $s^h_n$, and $p^h_n$   in terms of the respective basis functions at time level $n$.

\section{Newton-Krylov Solver} \label{sec:newton_krylov_solver}
The nonlinear system \eqref{eq:IBVP_G_res} is solved by means of Newton's method, with a backtracking strategy added for robustness.
The solution at time $t_n$ is computed as follows.
Given an initial guess $\mathbf{x}_n^0$, for $k = 0,1,\ldots$, until convergence
\begin{equation}
\begin{aligned}
 &\text{solve} &&A(\mathbf{x}_n^{(k)}) \Delta \mathbf{x} = - F(\mathbf{x}_n^{(k)}), \\
 &\text{set}   &&\mathbf{x}_n^{(k+1)} = \mathbf{x}_n^{(k)} + \lambda \Delta \mathbf{x},
\end{aligned}
\label{eq:newton}
\end{equation}
where $A(\mathbf{x}_n^{(k)}) = (\partial F/\partial \mathbf{x}_n)(\mathbf{x}_n^{(k)})$ is the Jacobian matrix associated with the nonlinear residual function $F$, and $\lambda \in (0,1]$ is an appropriately chosen line-search parameter.
For convenience of notation, we omit from now on to specify that $A$ is evaluated at $\mathbf{x}_n^{(k)}$.
Clearly, at each nonlinear iteration $k$, the solution of a linear system with $A$ is required.

The linearization of \cref{eq:IBVP_G_res} produces a Jacobian system with an inherent $3 \times 3$ block structure
\begin{align}
A = \begin{pmatrix}
A_{uu} & A_{us} & A_{up} \\
A_{su} & A_{ss} & A_{sp} \\
A_{pu} & A_{ps} & A_{pp}
\end{pmatrix}. \label{eq:block_linear_system}
\end{align}
This system has size proportional to the number of vertices (three displacement degrees of freedom per vertex) and cells (one saturation and one pressure degree of freedom per cell) in the computational mesh.
For detailed expressions of the sub-blocks in $A$, we refer the reader to \cite{White18}. 
Briefly, we emphasize the properties of the three diagonal blocks that motivate choices in designing the preconditioning operator described in \cref{sec:mgr}.
Specifically:
\begin{itemize}
\item $A_{uu}$ is the elasticity block and has the structure of a discrete elliptic operator;
\item $A_{ss}$ is the saturation block that, in the abscence of capillarity effects, has the structure of a discrete time-dependent hyperbolic problem;
\item $A_{pp}$ is the pressure block that,  similar to the elastic block, has the structure of a discrete elliptic operator.
\end{itemize}  
In this work, the linear system with matrix $A$ is solved iteratively with generalized minimal residual (GMRES) \cite{Saad86}, a Krylov subspace method designed for nonsymmetric systems.
Since Krylov methods' practical convergence depends 
on the availability of an effective preconditioner, we introduce the preconditioning operator $M$ and replace the linear system in \cref{eq:newton} with the right preconditioned system,
\begin{align}
AM^{-1} \Delta \mathbf{y} = -F(\mathbf{x}_n^{(k)}),
\end{align}
where $\Delta \mathbf{x} = M^{-1} \Delta \mathbf{y}$.
In the following section, we describe an algebraic method to construct $M$ given a matrix $A$ with the structure specified in \cref{eq:block_linear_system}.

\section{Multigrid Reduction}\label{sec:mgr}
The idea of MGR has been around for a long time, tracing back to the work of Ries and Trottenberg \cite{Ries79,Ries83}. 
Recently, it has gained more attention through the work on multigrid reduction in time by Falgout et al. \cite{Falgout14,Falgout16}. 
MGR has also been applied successfully for problems in reservoir simulation and multiphase flow in porous media with phase transitions \cite{Bui18,Wang17}.
A major advantage of the MGR approach is that it is an algebraic method and unlike geometric multigrid, it can be applied to general geometries and grid types.
In this section, we first summarize the approach for the case of two-level reduction and then present the general multi-level reduction algorithm.

\subsection{Two-grid Reduction Scheme}
For a matrix $A$ of size $N\times N$, we define a partition of the row indices of the matrix into C-points and F-points. 
The C-points play a role analogous to the points on a coarse grid, and the F-points belong to the set that is the complement of the C-points.
It is important to note that this partitioning is different from the one normally used in standard multigrid methods, in which the F-points correspond to all points on the fine grid, i.e. the set of F-points contains the set of C-points.
In multigrid reduction, the C-points and F-points belong to non-overlapping sets.
Following \cite{Falgout14}, using such CF-splitting we have
\begin{align}
A = \begin{pmatrix}
A_{FF} & A_{FC} \\
A_{CF} & A_{CC}
\end{pmatrix} = \begin{pmatrix}
I_{FF} &  0\\
A_{CF}A_{FF}^{-1} & I_{CC}
\end{pmatrix} \begin{pmatrix}
A_{FF} & 0 \\
0 & S
\end{pmatrix} \begin{pmatrix}
I_{FF} & A_{FF}^{-1}A_{FC}  \\
0& I_{CC}
\end{pmatrix},
\end{align}
where $I_{CC}$ and $I_{FF}$ are identity matrices and $S = A_{CC} - A_{CF} A_{FF}^{-1} A_{FC}$ is the Schur complement.
We can define the ideal interpolation and restriction operators by 
\begin{align}
P = \begin{pmatrix}
-A_{FF}^{-1} A_{FC}\\
I_{CC}
\end{pmatrix}, \hspace{5mm} R = \begin{pmatrix}
-A_{CF}A_{FF}^{-1} &I_{CC}
\end{pmatrix}.
\end{align}
Additionally, define the injection operator as $Q = \begin{pmatrix} I_{FF} \\ 0 \end{pmatrix}$.
Then since $A_{FF} = Q^TAQ$ and $S = RAP$, it is simple to derive that
\begin{equation}
A^{-1}=P(RAP)^{-1}R+Q(Q^TAQ)^{-1}Q^T,
\end{equation}
and 
\begin{align}
  0 = I - A^{-1}A
  &= I - P(RAP)^{-1}RA-Q(Q^TAQ)^{-1}Q^TA\label{eq:MGR-add}\\
  &=(I - P(RAP)^{-1}RA)(I-Q(Q^TAQ)^{-1}Q^TA)\label{eq:MGR-mul1}\\
  &=(I-Q(Q^TAQ)^{-1}Q^TA)(I - P(RAP)^{-1}RA),\label{eq:MGR-mul2}
\end{align}
where the equivalence occurs since $RAQ=Q^TAP=0$.
This identity defines the two-level multigrid method with the ideal Petrov-Galerkin coarse-grid correction $(RAP)^{-1}$ and the F-relaxation $Q(Q^TAQ)^{-1}Q^T$: (i) \Cref{eq:MGR-add} is the additive MGR identity and (ii) \cref{eq:MGR-mul1,eq:MGR-mul2} are multiplicative identities with pre-smoothing and post-smoothing F-relaxation, respectively.
However, constructing ideal interpolation and restriction operators is impractical. Similarly, computing the coarse-grid correction exactly is expensive, so we need to approximate these operators.
In practice, MGR methods use a scalable solver such as AMG for the coarse-grid solve, and replace the ideal restriction and prolongation $R$ and $P$ with
\begin{equation}
\label{eq:RPoperator1}
\tilde{P}= 
\begin{pmatrix}
W_{p}
\\ I_{CC},
\end{pmatrix},\quad
\tilde{R} =
\begin{pmatrix}
W_{r} & I_{CC}
\end{pmatrix}.
\end{equation}
where
\begin{align}
W_{r} \approx -A_{CF}A_{FF}^{-1}, \hskip2ex W_{p} \approx - A_{FF}^{-1} A_{FC}.
\end{align}
There are many ways to construct these approximations.
One simple choice is to use an injection operator for restriction and a Jacobi approach for interpolation
\begin{align}
W_r = 0, \hskip2ex W_p = -D_{FF}^{-1} A_{FC}, \label{eq:injective_rp}
\end{align}
where $D_{FF} = \text{diag}(A_{FF})$.
Then the coarse grid operator $A_{h} = \tilde{R}A\tilde{P}$ can also be considered as an approximation to the Schur complement $S$.
Besides the choices in \cref{eq:injective_rp}, one can also choose to use Jacobi approach for restriction, that is $W_r = - A_{CF} D_{FF}^{-1}$.
Another option is to construct $A_{FF}^{-1}$ using incomplete factorizations (ILU) or sparse approximate inverse techniques, such as sparse approximate inverse (SPAI) \cite{Grote97}, factored sparse approximate inverse (FSAI) \cite{Ferronato14}, or minimal residual (MR) \cite{Chow98}.
Although these methods could provide a better approximation to $A_{FF}^{-1}$, and therefore better approximations for the restriction and interpolation operators, they tend to make these operators dense.
The resulting coarse grid also becomes dense and unamenable to AMG.
One can certainly apply a dropping strategy to keep such $\tilde{P}$ and $\tilde{R}$ sparse, but in practice, the potential improvement in performance using approximate inverse methods is usually offset by the cost to construct the approximation, which makes simple methods like Jacobi more appealing.

In general, we define the MGR operator with either pre-smoothing or post-smoothing F-relaxation by
\begin{align}
  I - M_{MGR}^{-1}A &= (I - \tilde{P}M_{CC}^{-1}\tilde{R}A)(I-M_{FF}^{-1}A) , \label{eq:MGR-pre}\\
  I - M_{MGR}^{-1}A &= (I-M_{FF}^{-1}A)(I - \tilde{P}M_{CC}^{-1}\tilde{R}A) , \label{eq:MGR-post}
\end{align}
where $M_{CC} = (\tilde{R}A\tilde{P})$ is the coarse-grid correction and $M_{FF}^{-1}$ is the F-relaxation smoo\-th\-er. 
Additionally, similar to AMG methods, one can also apply a global smoothing step that extends to all the unknowns, not just the F-points.
For the global smoother $M_{glo}^{-1}$, various methods including (block) Jacobi, (block) Gauss-Seidel, or ILU, can be used.
The inclusion of this step can help eliminate error modes that both the F-relaxation and coarse-grid correction may have missed.
The application of the two-grid MGR scheme consisting of a global smoother and an F-relaxation followed by a coarse-grid correction can be summarized as shown in \cref{algo:two_level_mgr_NC}.


\begin{algorithm}
  \caption{Two-grid MGR preconditioner with presmoothing, $z = M_{MGR}^{-1} v$.}\label{algo:two_level_mgr_NC}
  \begin{algorithmic}[1]
  \Function{\tt {MGR}}{$A, v$}
    \State $z = M_{glo}^{-1} v $ \Comment{Global Relaxation}
    \State $z \leftarrow z + QM_{FF}^{-1}Q^{T} (v - Az)$ \Comment{F-Relaxation}
    \State $r_C = \tilde{R}(v - Az)$ \Comment{Restrict residual}
    \State $M_{CC} e_C = r_C$ \Comment{Solve coarse-grid error problem with AMG}
    \State $e = \tilde{P} e_C$ \Comment{Interpolate coarse error approximation}
    \State $z \leftarrow z + e$ \Comment{Apply correction}
    \State \Return{$z$}
  \EndFunction
  \end{algorithmic}
\end{algorithm}


Balancing the quality of the approximation to the Schur-complement and the convergence of the coarse-grid solve is key to the success of MGR.
One extreme is to design a coarse grid that is perfectly suitable for AMG.
Assuming, for example, that the block $A_{CC}$ comes from a scalar elliptic PDE and $A_{CC}$ is SPD, then one can choose $W_p = W_r = 0$ and the coarse grid becomes $RAP = A_{CC}$.
In this case, the convergence of the coarse grid solve is optimal, but the approximation of the Schur-complement far from ideal, since the coarse grid neither takes into account any information from the F-points nor the coupling between the C and F points.
At the other extreme, one can use the exact Schur-complement as the coarse grid by choosing $W_r = -A_{CF}A_{FF}^{-1}$ and $W_p = 0$.
However, because of the exact inversion of $A_{FF}^{-1}$, the coarse grid is dense.
Furthermore, since the F-points and C-points actually represent equations obtained from the discretization of different continuous physical models, capturing the coupling between them on the coarse grid can lead to loss of ellipticity, which can make the coarse-grid solve with AMG ineffective.
Thus, finding a good approximation of the Schur-complement that is still amenable to AMG methods is essential.

\begin{remark}
The appeal of the MGR approach is that it provides a general framework for choosing the coarse/fine grids, the interpolation and restriction operators, and the solvers for the F-relaxation and coarse-grid correction.
As an example, it was shown in \cite{Bui18,Wang17} that one can recast any CPR-AMG strategy \cite{Cao05,Gries14,Lacroix03,Liu15,Scheichl03,Stueben07,Zhou12} or block preconditioner \cite{Bui17} used in reservoir simulation as a particular variant of the two-grid MGR reduction scheme by appropriately defining the different components of the algorithm, namely prolongation, restriction and smoothing operators.
\end{remark}

\subsection{A general multi-level MGR algorithm} 
One can replace the coarse grid solve in \cref{algo:two_level_mgr_NC} with a two-level MGR scheme and apply the method recursively to obtain a multi-level MGR algorithm.
The general application of the MGR V-cycle with global smoothing is summarized in \cref{algo:multi_level_mgr_NC}, where the hierarchy of coarse grid operators, i.e. $A_{l+1} = \tilde{R}_{l} A_{l} \tilde{P}_{l}$, is assumed to be computed for each level $l$.


\begin{algorithm}
  \caption{General multi-level MGR preconditioner, $z = M_{l,MGR}^{-1} v$.}\label{algo:multi_level_mgr_NC}
  \begin{algorithmic}[1]
  \Function{\tt {MGR}}{$A_{l}, v_{l}$}
    \If{$l$ is the coarsest level}
      \State $A_{l} z_{l} = v_{l}$ \Comment{Solve coarse-grid error problem with AMG}
    \Else
      \State $z_{l} = M_{l,glo}^{-1} v_{l} $ \Comment{Global Relaxation}
      \State $z_{l} \leftarrow z_{l} + Q_{l}M_{l,FF}^{-1}Q_{l}^{T} (v_{l} - Az_{l})$ \Comment{F-Relaxation}
      \State $r_{l+1} = \tilde{R}_{l} (v_{l} - A_{l}z_{l})$ \Comment{Restrict residual}
      \State  $e_{l+1} = {\tt {MGR}}(A_{l+1}, r_{l+1})$  \Comment{Recursion} 
      \State $e_{l} = \tilde{P}_{l} e_{l+1}$ \Comment{Interpolate coarse error approximation}
      \State $z_{l} \leftarrow z_{l} + e_{l}$ \Comment{Apply correction}
   	\EndIf
   	\State \Return{$z_{l}$}
  \EndFunction
  \end{algorithmic}
\end{algorithm}

Based on \cref{algo:multi_level_mgr_NC}, W- and F-cycle versions of the MGR algorithm can also be defined \cite{Ries83}. 
Note that the Schur-complement $S$ is approximated by the triple product $RAP = A_{CC} - A_{CF}D_{FF}^{-1}A_{FC}$ in the classical two-grid reduction scheme in \cref{algo:two_level_mgr_NC}.
Even though we have introduced a sparse approximation to $S$ by replacing $A_{FF}$ with its diagonal $D_{FF}$, i.e. $A_{FF}^{-1} \approx D_{FF}^{-1}$, in a multi-level reduction scheme, the coarse grid can still become dense or unsuitable for standard AMG because the correction term $A_{cor} = A_{CF}D_{FF}^{-1}A_{FC}$ involves a matrix-matrix product. 

In this work, we develop a dropping strategy for $A_{cor}$ to keep the coarse grid sparse as well as suitable for AMG. One approach is to drop all entries of $A_{cor}$ that are smaller than a prescribed tolerance. Here, we use a different strategy based on a maximum number of non-zero values per row. Specifically, we choose to keep only $N_{\text{max}}$ entries with largest absolute values on each row.
%
%
To preserve at least some information of the first level of reduction, however, we always keep the diagonals of the sub-blocks in $A_{cor}$.
For instance, in a three-level reduction scheme, in the first-level reduction, $A_{cor}$ has $2\times 2$ block structure, and applying maximum dropping (i.e. using an extremely large tolerance or $N_{\text{max}} = 0$), $A_{cor}$ is still a $2 \times 2$ block matrix, whose sub-blocks are diagonal matrices.
Using this dropping strategy results in a non-Galerkin coarse grid
\begin{align}
S = A_{CC} - \mathbf{G}(A_{CF}D_{FF}^{-1}A_{FC}).
\label{eq:drop_operator}
\end{align}
where $\mathbf{G}$ is a sparsifying operator that performs one of the aforementioned dropping strategies.
So far we only assume that a CF-splitting of the rows is given.
How to choose such a splitting is dependent on the problem and it is up to the user to make the decision.
However, as a general principle, it is usually a good idea to choose 
a CF-splitting so that the final coarse grid corresponds to the variable associated with an elliptic equation, e.g. pressure, since we want to solve the coarse grid using an efficient method such as standard AMG.
In the next section, we show how to choose an appropriate CF-splitting at each level of reduction for our multiphase poromechanical problem.

\subsection{MGR for Multiphase Poromechanics}
We propose a three-level MGR reduction scheme to precondition the Jacobian matrix \eqref{eq:block_linear_system}.
For the first level of reduction, we aim at decoupling the mechanics sub-problem from the flow.
Therefore, we assign all the displacement unknowns as F-points while both saturation and pressure unknowns are labeled as C-points.
This leads to the following partitioning
\begin{align}
A &= \left(\begin{array}{c|cc}
A_{uu} & A_{us} & A_{up} \\
\hline
A_{su} & A_{ss} & A_{sp} \\
A_{pu} & A_{ps} & A_{pp}
\end{array} \right) \begin{array}{cc}
& F\\
& C\\
& C
\end{array}.
\end{align}
Then $A_{FF} \equiv A_{uu}$ and the coarse grid $A_{CC}$, which corresponds to the flow sub-problem, has the $2\times 2$ block structure
\begin{align}
A_{CC} = \begin{pmatrix}
A_{ss} & A_{sp} \\
A_{ps} & A_{pp}
\end{pmatrix}.
\end{align}
For the F-relaxation step, we need to solve the elasticity problem involving the elliptic operator $A_{uu}$.
Here, we use one AMG V-cycle. 
Because of the vectorial nature of the elasticity operator, this is the most expensive part of the setup phase.
Also, given the global system size, we ignore the first-level global relaxation step.
Using the interpolation and restriction operators specified in \cref{eq:injective_rp} combined with the dropping strategy defined in \Cref{eq:drop_operator} yields the following first level coarse grid
\begin{equation}
\begin{aligned}
S_1 &=  \begin{pmatrix}
A_{ss} & A_{sp} \\
A_{ps} & A_{pp}
\end{pmatrix} - 
\mathbf{G}
\left(
\begin{pmatrix}
A_{su} \\
A_{pu}
\end{pmatrix} D_{uu}^{-1} 
\begin{pmatrix}
A_{us} & A_{up}
\end{pmatrix} \right) \\
&= \begin{pmatrix}
\tilde{A}_{ss} & \tilde{A}_{sp} \\
\tilde{A}_{ps} & \tilde{A}_{pp}
\end{pmatrix}. \label{eq:coarse_grid_flow_sparse}
\end{aligned}
\end{equation}
For our multiphase poromechanics problem, we use $N_{\text{max}} = 4$ for $\mathbf{G}$.
Again, we emphasize the flexibility of our framework as it allows for experimenting with different choices of $\mathbf{G}$.
For example, choosing an appropriate $\mathbf{G}$, we can mimic the fixed-stress preconditioner developed in \cite{White18}.
\par The second reduction step is essentially a CPR approach that is embedded within a multigrid reduction framework.
Hence, we label saturation unknowns as F-points and pressure unknowns as C-points in $S_1$: 
\begin{align}
\renewcommand*{\arraystretch}{1.2}
S_1 &= \left(\renewcommand*{\arraystretch}{1.25}\begin{array}{c|c}
\tilde{A}_{ss} & \tilde{A}_{sp} \\
\hline
\tilde{A}_{ps} & \tilde{A}_{pp}
\end{array} \right) \begin{array}{cc}
& F\\
& C
\end{array}.
\end{align}
Again, using interpolation and restriction operators in \cref{eq:injective_rp}, we obtain the second-level Schur-complement 
\begin{align}
S_2 = \tilde{A}_{pp} - \tilde{A}_{ps} D_{ss}^{-1} \tilde{A}_{sp}, \label{eq:coarse_grid_pressure}
\end{align}
where $D_{ss} = \text{diag} (\tilde{A}_{ss})$.
In other physics-based approaches commonly used in reservoir simulation, one can seek to further sparsify the Schur-complement. 
For example, in a \textit{Quasi-IMPES reduction} scheme, the block $\tilde{A}_{ps}$ is also replaced by its diagonal $D_{ps} = \text{diag} (\tilde{A}_{ps})$.
This approximation ensures that the matrix sparsity pattern of $A_{pp}$ coming from the original finite volume stencil is preserved and the resulting Schur-complement is near elliptic.
In the MGR approach, however, no further sparse approximation is needed for this level since the flow part is relatively small compared to the elasticity block and the coarse grid generated in \cref{eq:coarse_grid_pressure} is still well-suited for AMG.
At the second level, the F-relaxation involving the $\tilde{A}_{ss}$ block is done using a simple Jacobi relaxation.
However, the second-level global smoothing step is required to reduce the error associated with the hyperbolic component of the flow subproblem.
Indeed, the global smoothing plays a key role particularly in the later stage of the simulation when the pressure field approaches steady-state conditions and the multiphase flow and transport process transitions to an advection dominated regime.
The need for a robust global smoother will become clear through numerical results presented in the next section.

\begin{remark}
Even though we formally present the multigrid reduction framework for the $3\times 3$ system in field-ordered form, in our implementation, 
the input matrix has interleaved ordering for saturation and pressure.
This choice produces a sparsity pattern in which dense 
$2\times 2$ blocks appear for the first-level coarse grid in \cref{eq:coarse_grid_flow_sparse}.
%
%
Block versions of relaxation or incomplete factorization preconditioners are therefore appealing, as dense multiplication and inversion operations can be applied to the small blocks.
\end{remark}

\begin{remark}
Common strategies for the second-level global smoothing step include block relaxation methods (e.g. Jacobi, Gauss-Seidel) or incomplete factorizations (e.g. ILU(k), ILUT).
In this work, we explore two options.
The first option uses several sweeps of hybrid block Gauss-Seidel (HBGS).
The second option uses one sweep of processor-local, pointwise ILU(k) \cite{Chow15}.
\end{remark}

\section{Numerical Results}\label{sec:numerical_results}
We perform numerical experiments to test the performance of the MGR preconditioner on two problems: 
(1) a weak scaling study for a simple synthetic configuration; and (2) a strong scaling study using a realistic, 
highly heterogeneous reservoir based on the SPE10 \cite{Christie01} example.  {Both examples have been designed as community benchmark problems and exhibit tight-coupling between displacement, pressure, and saturation fields.  Problem specifications are described in detail in \cite{White18} and so are only briefly reported below.

In this study, the simulator is 
provided by \textit{Geocentric}, which utilizes the \textit{deal.ii} Finite Element Library \cite{Bangerth07} for 
discretization functionality. It also provides a direct interface with MGR, which is implemented as a separate 
solver in \textit{hypre} \cite{Falgout02}. All the numerical experiments were run on \textit{Quartz}, a cluster 
at the Lawrence Livermore Computing Center with 1344 nodes containing two Intel Xeon E5-2695 18-core 
processors sharing 128 GiB of memory on each node with Intel Omni-Path interconnects between nodes. 
We use pure MPI-based parallelism.

For the elastic block $A_{uu}$, we use one V-cycle of BoomerAMG \cite{Henson00}, with an unknown approach for a system of three PDEs, with one level of 
aggressive coarsening, one sweep of hybrid forward $l_1$-Gauss-Seidel \cite{Baker11} for the down cycle and 
one sweep of hybrid backward $l_1$-Gauss-Seidel for the up cycle. The coarsest grid is solved directly with 
Gaussian elimination. The MGR coarse-solve in \cref{eq:coarse_grid_pressure} also uses BoomerAMG with the same smoother configuration, but for a 
scalar problem and a Hybrid Modified Independent Set (HMIS) coarsening strategy \cite{DeSterck06}. For the 
global smoother, we use one step of processor-local, pointwise ILU(1).

\subsection{Staircase Benchmark}

\begin{figure}
  \centering
  \begin{subfigure}[t]{0.42\textwidth}
    \centering
    \includegraphics[width=\linewidth]{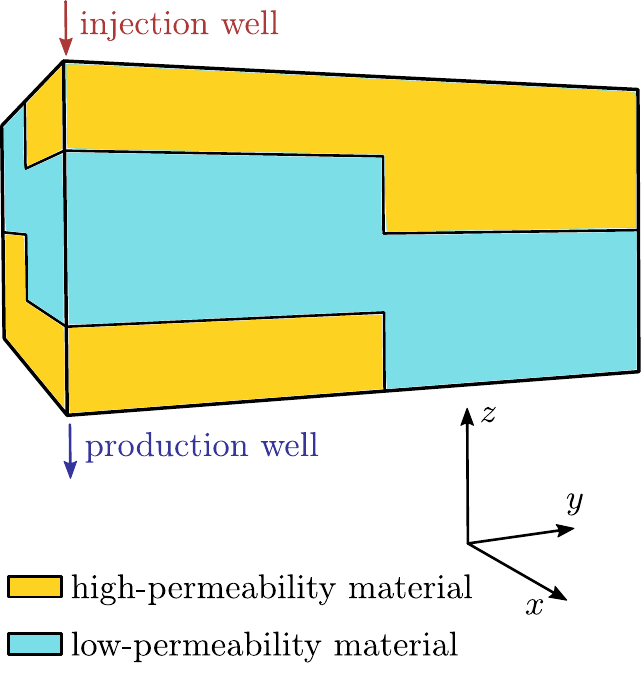}
    \caption{Sketch of the simulated domain.}
    \label{fig:staircase_a}
  \end{subfigure}
  \hfill
  \begin{subfigure}[t]{0.42\textwidth}
    \centering
    \includegraphics[width=\linewidth]{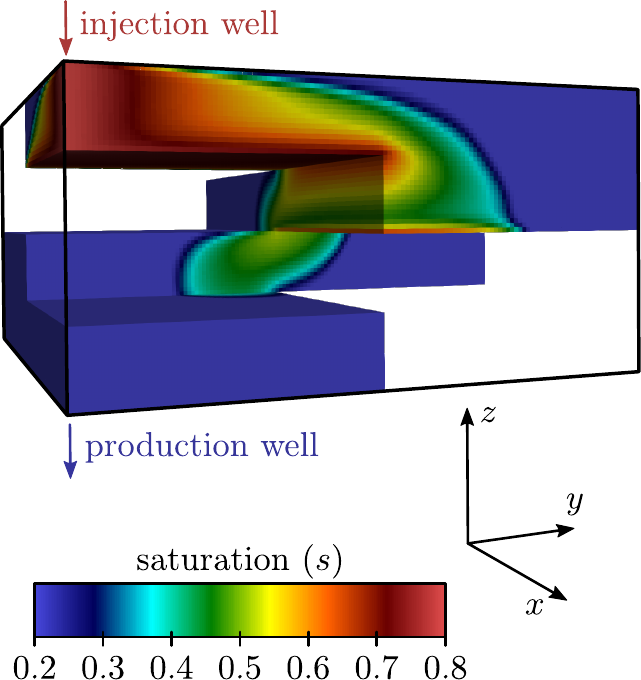}
    \caption{Saturation field.}
    \label{fig:staircase_b}
  \end{subfigure}
  \caption{Staircase benchmark, showing basic geometry and resulting saturation field within the high-permeability channel at $t=92$ days.}
  \label{fig:staircase}
\end{figure}

The configuration of the first test problem is illustrated in Figure~\ref{fig:staircase}. A highly-permeable channel winds its way in a ``staircase" fashion through a lower-permeability host rock.  A denser, wetting phase is injected through a well at the top corner, leading to a saturation plume driven by gravity and pressure that migrates towards a production well in the lower corner.  The whole system is deformable and exhibits significant poromechanical coupling.  Visualizations of the resulting pressure and deformation fields have been omitted for brevity.  A detailed specification of mesh geometry, material properties, and boundary conditions can be found in \cite{White18}.

\begin{table}[t]
\centering
\caption{Weak scaling performance for the staircase example.} \label{tab:weak_scaling}
\includegraphics[width=0.9\textwidth]{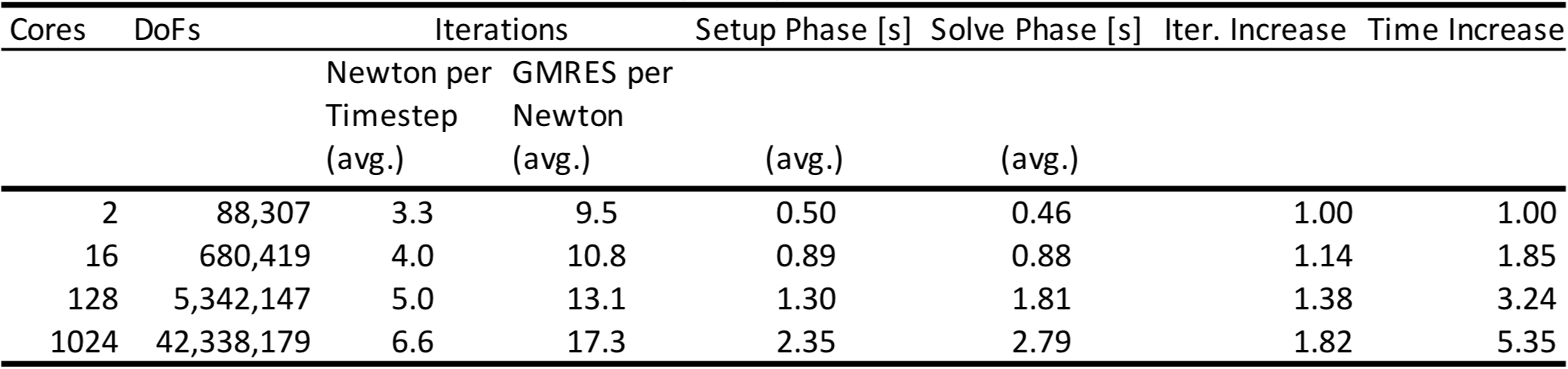}
\end{table}

\Cref{tab:weak_scaling} shows the results for a weak scaling study using the staircase example. We keep the 
number of degrees of freedom per core constant at 44k and increase the number of cores from 2 to 1024. The 
global problem size grows $8^3$ times from 88k to 42M. Due to the inherent nonlinearity, we observe an increase in 
the number of Newton iterations per time step as the mesh is refined. The average number of GMRES iterations 
per Newton step, however, only experiences a modest growth as desired. Even though we can use a more 
complex smoother in place of the hybrid $l_1$-Gauss-Seidel solves for the elasticity block and drive down 
the number of iterations, that will come at the expense of run-time performance. In general, we find that the 
$l_1$-Gauss-Seidel smoother strikes a good balance between iteration counts and run time. Similar to the 
number of iterations, the total run time, including both the setup and solve phases, also exhibits some growth, 
but again, the result is quite satisfactory even for large core counts. The increase in the run time can be attributed 
to communication costs in the MGR setup and solve phases, since the actual number of degrees of freedom per 
core is fairly small. Overall, however, the MGR framework provides a good platform for scalable performance.

\subsection{SPE10-based Benchmark}

\begin{figure}
  \centering
  \begin{subfigure}[t]{0.48\textwidth}
    \centering
    \includegraphics[width=\linewidth]{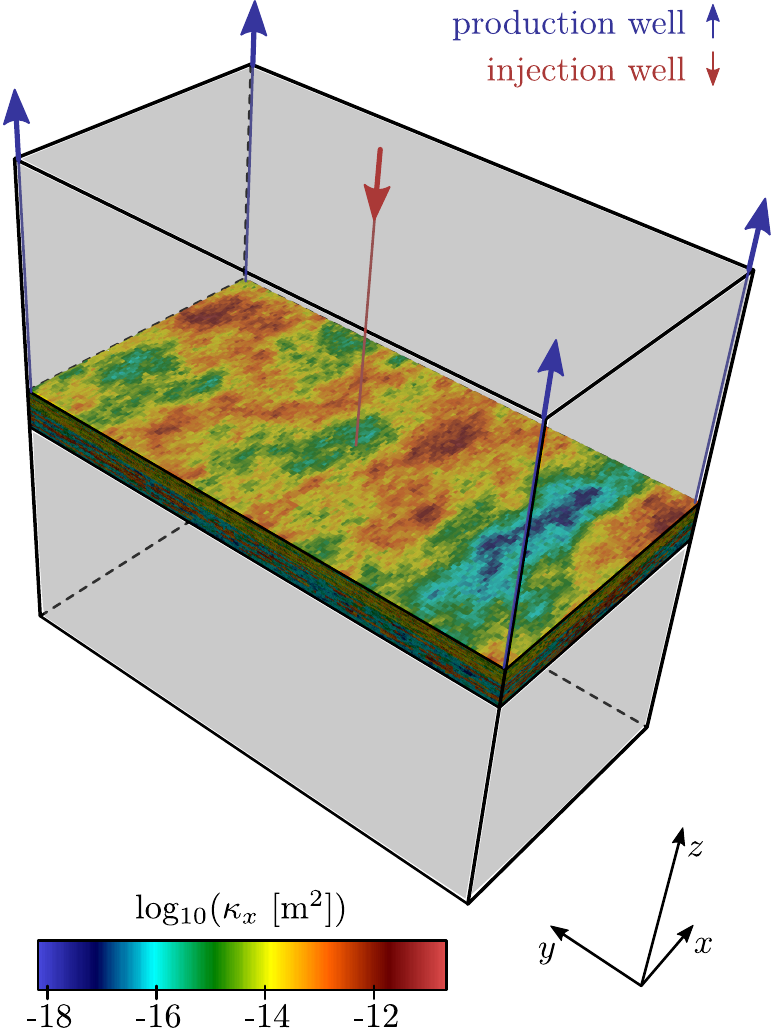}
    \caption{Sketch of the simulated domain with horizontal permeability field in the produced reservoir. }
    \label{fig:SPE10_a}
  \end{subfigure}
  \hfill
  \begin{subfigure}[t]{0.48\textwidth}
    \centering
    \includegraphics[width=\linewidth]{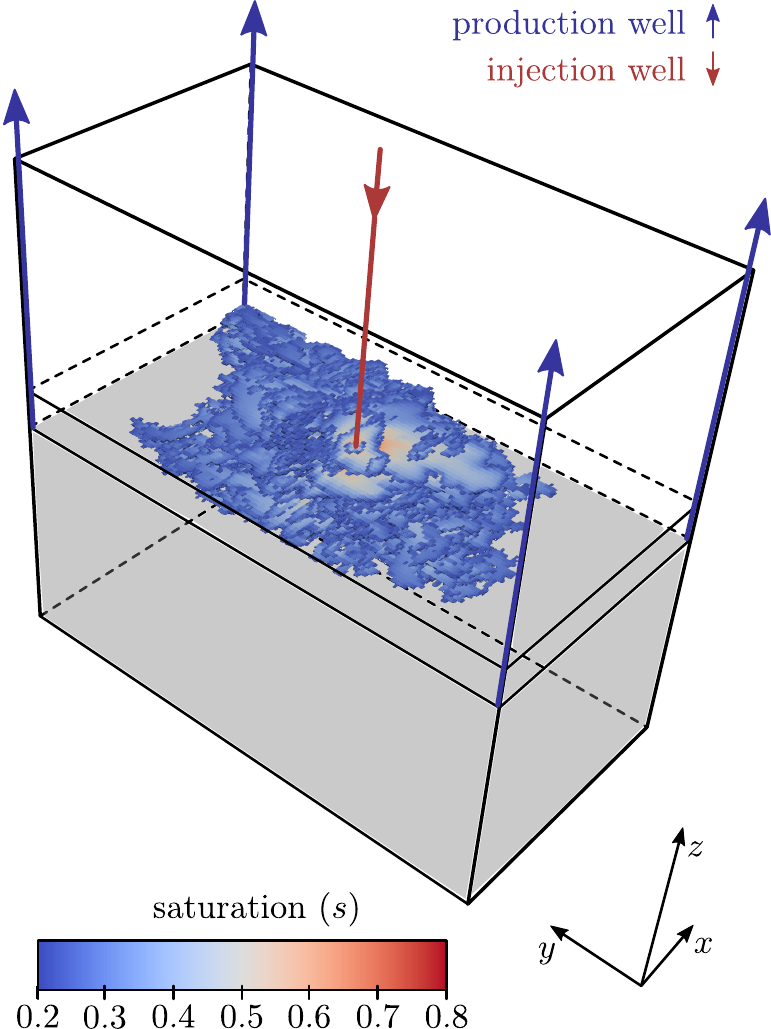}
    \caption{Saturation field at $t$=98 days.}
    \label{fig:SPE10_b}
  \end{subfigure}
  \caption{SPE10-based benchmark: The original SPE-10 reservoir is embedded in a larger poromechanical domain to provide realistic mechanical boundary conditions.}
  \label{fig:SPE10}
\end{figure}

\begin{table}[t]
\centering
\caption{Strong scaling performance for the SPE10-based problem.} \label{tab:strong_scaling}
\includegraphics[width=0.9\textwidth]{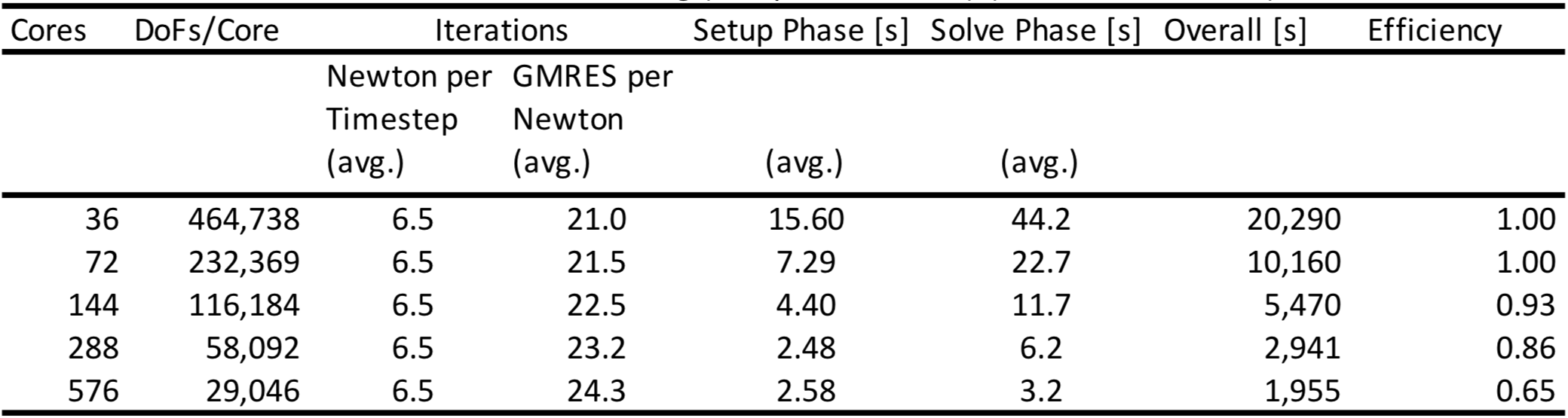}
\end{table}

We also perform a strong scaling study on a more realistic benchmark problem derived from the second 
model of the SPE10 Comparative Solution Project (Figure~\ref{fig:SPE10})  \cite{Christie01}.  The original SPE10 permeability and porosity fields are now treated as a poromechanical medium.  These geostatistically generated fields exhibit both severe heterogeneity and anisotropy.  In the current poromechanical benchmark, the reservoir itself is also embedded in a larger domain---with caprock and underburden---to provide more realistic boundary conditions.  Water is injected through a central well, while fluids are produced from four wells at the corners of the domain.  Mesh, material property, and boundary condition specifications are reported in \cite{White18}.  Note that the well control conditions differ from the original SPE10 model to avoid well impacts on the linear solver.  The treatment of well degrees-of-freedom within the linear solver is a critical issue, but is deliberately left out-of-scope for the current contribution.  We remark, however, that the MGR approach provides a flexible framework to treat this additional complexity. 

The resulting discrete problem has 16.7M degrees-of-freedom. 
We keep the problem size fixed and divide the work across an increasing number of 
compute cores. The results are shown in \cref{tab:strong_scaling}. Again, we observe only minor growth in the 
number of GMRES iterations with larger core counts. Similar to the weak scaling case, the reason for this growth 
is the use a hybrid $l_1$-Gauss-Seidel smoother in AMG solves for the elasticity block and the coarse grid. 
Good overall timing efficiency is also achieved up to 288 cores. For 576 cores, even though we still get good 
efficiency for the solve phase, there is a noticeable increase in the setup time because the problem size on each 
core becomes very small, i.e. about 17k total and less than 6k degrees of freedom for the elastic block and the coarse 
grid, respectively. Consequently, the majority of the time is spent in communication while not much computation 
is performed. However, the results still indicate that one can use the proposed framework with a large number of 
processors to efficiently reduce the long simulation time for challenging problems with highly heterogeneous media.

\subsection{Effect of global smoother}

\begin{figure}
\centering
\begin{subfigure}[t]{0.44\textwidth}
\centering
\includegraphics[width=\textwidth]{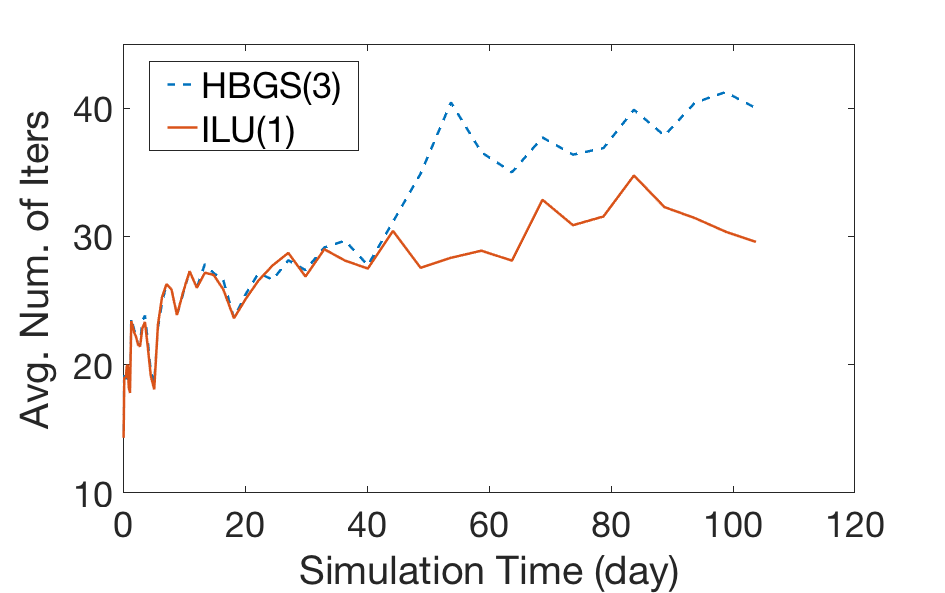}
\captionsetup{justification=centering}
\caption{Number of GMRES iterations.} \label{fig:smoother_comparison_iters}
\end{subfigure}
\begin{subfigure}[t]{0.45\textwidth}
\centering
\includegraphics[width=\textwidth]{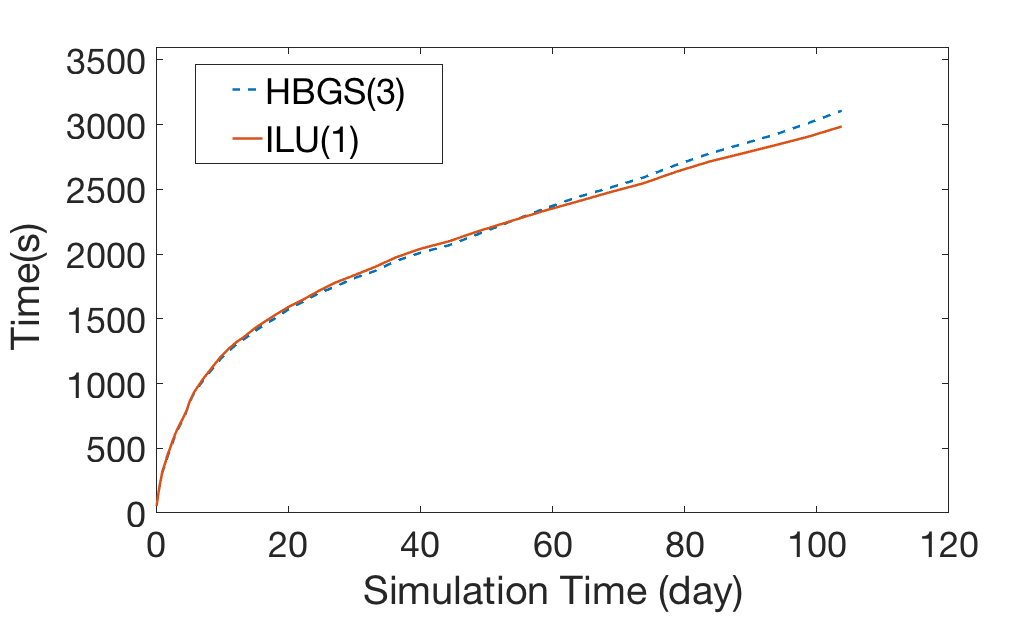}
\captionsetup{justification=centering}
\caption{Run time.}  \label{fig:smoother_comparison_timing}
\end{subfigure}
\captionsetup{justification=centering}
\caption{Effects of the global smoother at different time steps for the hybrid block Gauss-Seidel and processor-local ILU(1).}
\end{figure}

As we have mentioned earlier, the performance of MGR is dependent on the effectiveness of the solvers for each 
component of the algorithm. In general, changing the configuration for one component, e.g. smoother choices for 
the F-relaxation or coarsening strategies for the coarse-grid AMG solve, would result in a different 
number of GMRES iterations. However, the effect could also be quite subtle and not manifest itself until the 
underlying property of the problem changes. For multiphase poromechanics simulations, early times are typically dominated by elliptic effects associated with the pressure and displacement fields, while at late times the hyperbolic effects associated with the saturation field become significant. 

Here we explore the effectiveness of different global smoothers on 
the multiphase flow system as the simulation progresses. The first option uses three sweeps of HBGS, and the 
second option uses a single sweep of processor-local ILU(1). As one can see from \cref{fig:smoother_comparison_iters}, 
there is no apparent difference between the two smoothers until about 45 days of injection, when the number of 
iterations for the HBGS method increases sharply and continues to stay high. In contrast, ILU(1) is less sensitive. Even though the number of iterations also rises slightly around 85 days, it starts to decrease for the last 
period at the end of the simulation. We also plot the total time of HBGS(3) 
against ILU(1) in \cref{fig:smoother_comparison_timing}. It is clear that even though 
ILU(1) takes slightly more time in the beginning (mainly due to higher setup cost), the trade off is worthwhile 
thanks to its robustness, which leads to a modest reduction in total time for the whole simulation. This observation is confirmed by the widespread use of incomplete factorization smoothers in the reservoir simulation community.


\section{Conclusion}\label{sec:conclusion}
In this work, we have presented an algebraic framework based on multigrid reduction for solving the linear system that comes from discretizing 
and linearizing the conservation equations governing multiphase flow coupled with poromechanics. This framework is flexible and allows us to 
construct different preconditioners based on different choices for CF-splitting strategies, interpolation and restriction operators, as well as solvers 
and smoothers. We have also developed a dropping strategy for computing the reduction onto the coarse-grid within the MGR V-cycle that 
captures the coupling between mechanics and flow and reduces the operator complexity at the same time. This results in an algebraic preconditioner 
that is robust and scalable for realistic and large-scale simulation of multiphase poromechanics.

Regarding future work, a number of improvements to the MGR framework could be explored. For example, constructing good approximations to the 
ideal interpolation and restriction operators that have low complexity remains a significant challenge. Also, it is unclear how one can choose an optimal 
coarse grid that is representative of the fine-grid system and at the same time still amenable to AMG in a multi-level reduction setting. Thus, better 
strategies for computing the non-Galerkin coarse grid introduced in this work are needed to improve robustness of the framework. Lastly, since MGR 
is designed to accommodate a wide range of coupled systems, we are looking into extending the approach to solve problems with non-isothermal 
flow and fractured media.

\bibliographystyle{siamplain}
\bibliography{../master_all}

\begin{thebibliography}{10}

\bibitem{Adler18}
{\sc F.~J. Adler, J.~H. Gaspar, X.~Hu, C.~Rodrigo, and L.~T. Zikatanov}, {\em
  Robust block preconditioners for biot's model}, in Domain Decomposition
  Methods in Science and Engineering {XXIV}, P.~E. Bj{\o}rstad, S.~C. Brenner,
  L.~Halpern, H.~H. Kim, R.~Kornhuber, T.~Rahman, and O.~B. Widlund, eds.,
  Springer International Publishing, 2018,
  \url{https://doi.org/10.1007/978-3-319-93873-8},
  \url{https://doi.org/10.1007/978-3-319-93873-8}.

\bibitem{Aziz79}
{\sc K.~Aziz and A.~Settari}, {\em Petroleum Reservoir Simulation}, Applied
  Science Publishers, 1979.

\bibitem{Baker11}
{\sc A.~H. Baker, R.~D. Falgout, T.~V. Kolev, and U.~M. Yang}, {\em Multigrid
  smoothers for ultraparallel computing}, {SIAM} Journal on Scientific
  Computing, 33 (2011), pp.~2864--2887,
  \url{https://doi.org/10.1137/100798806},
  \url{https://doi.org/10.1137/100798806}.

\bibitem{Bangerth07}
{\sc W.~Bangerth, R.~Hartmann, and G.~Kanschat}, {\em deal.{II}---a
  general-purpose object-oriented finite element library}, {ACM} Transactions
  on Mathematical Software, 33 (2007), pp.~24--es,
  \url{https://doi.org/10.1145/1268776.1268779},
  \url{https://doi.org/10.1145/1268776.1268779}.

\bibitem{Bergamaschi07}
{\sc L.~Bergamaschi, M.~Ferronato, and G.~Gambolati}, {\em Novel
  preconditioners for the iterative solution to {FE}-discretized coupled
  consolidation equations}, Computer Methods in Applied Mechanics and
  Engineering, 196 (2007), pp.~2647--2656,
  \url{https://doi.org/10.1016/j.cma.2007.01.013},
  \url{https://doi.org/10.1016/j.cma.2007.01.013}.

\bibitem{Bergamaschi12}
{\sc L.~Bergamaschi and {\'{A}}.~Mart{\'{\i}}nez}, {\em {RMCP}: Relaxed mixed
  constraint preconditioners for saddle point linear systems arising in
  geomechanics}, Computer Methods in Applied Mechanics and Engineering, 221-222
  (2012), pp.~54--62, \url{https://doi.org/10.1016/j.cma.2012.02.004},
  \url{https://doi.org/10.1016/j.cma.2012.02.004}.

\bibitem{Bui17}
{\sc Q.~M. {Bui}, H.~C. {Elman}, and J.~D. {Moulton}}, {\em Algebraic multigrid
  preconditioners for multiphase flow in porous media}, {SIAM} Journal on
  Scientific Computing, 39 (2017), pp.~{S662--S680},
  \url{https://doi.org/10.1137/16M1082652},
  \url{http://dx.doi.org/10.1137/16M1082652},
  \url{https://arxiv.org/abs/http://dx.doi.org/10.1137/16M1082652}.

\bibitem{Bui18}
{\sc Q.~M. {Bui}, L.~{Wang}, and D.~{Osei-Kuffuor}}, {\em Algebraic multigrid
  preconditioners for two-phase flow in porous media with phase transitions},
  Advances in Water Resources, 114 (2018), pp.~19--28,
  \url{https://doi.org/10.1016/j.advwatres.2018.01.027}.

\bibitem{Cao05}
{\sc H.~Cao, H.~A. Tchelepi, J.~H. Wallis, and H.~E. Yardumian}, {\em Parallel
  scalable unstructured {CPR}-type linear solver for reservoir simulation}, in
  {SPE} Annual Technical Conference and Exhibition, Society of Petroleum
  Engineers ({SPE}), 2005, \url{https://doi.org/10.2118/96809-ms},
  \url{http://dx.doi.org/10.2118/96809-MS}.

\bibitem{Chow15}
{\sc E.~Chow and A.~Patel}, {\em Fine-grained parallel incomplete {LU}
  factorization}, {SIAM} Journal on Scientific Computing, 37 (2015),
  pp.~C169--C193, \url{https://doi.org/10.1137/140968896},
  \url{https://doi.org/10.1137/140968896}.

\bibitem{Chow98}
{\sc E.~Chow and Y.~Saad}, {\em Approximate inverse preconditioners via
  sparse-sparse iterations}, {SIAM} Journal on Scientific Computing, 19 (1998),
  pp.~995--1023, \url{https://doi.org/10.1137/s1064827594270415},
  \url{https://doi.org/10.1137/s1064827594270415}.

\bibitem{Christie01}
{\sc M.~A. Christie and M.~J. Blunt}, {\em Tenth {SPE} comparative solution
  project: A comparison of upscaling techniques}, in {SPE} Reservoir Simulation
  Symposium, Society of Petroleum Engineers ({SPE}), 2001,
  \url{https://doi.org/10.2118/66599-ms},
  \url{http://dx.doi.org/10.2118/66599-MS}.

\bibitem{Cou04}
{\sc O.~Coussy}, {\em {Poromechanics}}, Wiley, Chichester, UK, 2004.

\bibitem{DeSterck06}
{\sc H.~{De Sterck}, U.~M. Yang, and J.~J. Heys}, {\em Reducing complexity in
  parallel algebraic multigrid preconditioners}, {SIAM} Journal on Matrix
  Analysis and Applications, 27 (2006), pp.~1019--1039,
  \url{https://doi.org/10.1137/040615729},
  \url{https://doi.org/10.1137/040615729}.

\bibitem{Falgout14}
{\sc R.~D. Falgout, S.~Friedhoff, T.~V. Kolev, S.~P. MacLachlan, and J.~B.
  Schroder}, {\em Parallel time integration with multigrid}, {SIAM} Journal on
  Scientific Computing, 36 (2014), pp.~C635--C661,
  \url{https://doi.org/10.1137/130944230},
  \url{http://dx.doi.org/10.1137/130944230},
  \url{https://arxiv.org/abs/http://dx.doi.org/10.1137/130944230}.

\bibitem{Falgout16}
{\sc R.~D. Falgout, T.~A. Manteuffel, B.~O'Neill, and J.~B. Schroder}, {\em
  Multigrid reduction in time for nonlinear parabolic problems}, tech. report,
  Lawrence Livermore National Laboratory, jan 2016,
  \url{https://doi.org/10.2172/1236132}, \url{https://doi.org/10.2172/1236132}.

\bibitem{Falgout02}
{\sc R.~D. Falgout and U.~Yang}, {\em {HYPRE}: a library of high performance
  preconditioners}, in Preconditioners, Lecture Notes in Computer Science,
  2002, pp.~632--641.

\bibitem{Ferronato14}
{\sc M.~Ferronato, C.~Janna, and G.~Pini}, {\em A generalized block {FSAI}
  preconditioner for nonsymmetric linear systems}, Journal of Computational and
  Applied Mathematics, 256 (2014), pp.~230 -- 241,
  \url{https://doi.org/10.1016/j.cam.2013.07.049},
  \url{http://www.sciencedirect.com/science/article/pii/S0377042713003944}.

\bibitem{GarKarTch16}
{\sc T.~A. Garipov, M.~Karimi-Fard, and H.~A. Tchelepi}, {\em {Discrete
  fracture model for coupled flow and geomechanics}}, Comput. Geosci., 20
  (2016), pp.~149--160, \url{https://doi.org/10.1007/s10596-015-9554-z}.

\bibitem{Gries14}
{\sc S.~Gries, K.~St\"{u}ben, G.~L. Brown, D.~Chen, and D.~A. Collins}, {\em
  Preconditioning for efficiently applying algebraic multigrid in fully
  implicit reservoir simulations}, {SPE} Journal, 19 (2014), pp.~726--736,
  \url{https://doi.org/10.2118/163608-pa},
  \url{https://doi.org/10.2118/163608-pa}.

\bibitem{Grote97}
{\sc M.~J. Grote and T.~Huckle}, {\em Parallel preconditioning with sparse
  approximate inverses}, {SIAM} Journal on Scientific Computing, 18 (1997),
  pp.~838--853, \url{https://doi.org/10.1137/s1064827594276552},
  \url{https://doi.org/10.1137/s1064827594276552}.

\bibitem{HagOsnLan12b}
{\sc J.~B. Haga, H.~Osnes, and H.~P. Langtangen}, {\em {On the causes of
  pressure oscillations in low-permeable and low-compressible porous media}},
  Int. J. Numer. Anal. Methods Geomech., 36 (2012), pp.~1507--1522,
  \url{https://doi.org/10.1002/nag.1062}.

\bibitem{Haga12}
{\sc J.~B. Haga, H.~Osnes, and H.~P. Langtangen}, {\em A parallel block
  preconditioner for large-scale poroelasticity with highly heterogeneous
  material parameters}, Computational Geosciences, 16 (2012), pp.~723--734,
  \url{https://doi.org/10.1007/s10596-012-9284-4},
  \url{https://doi.org/10.1007/s10596-012-9284-4}.

\bibitem{Henson00}
{\sc V.~E. Henson and U.~M. Yang}, {\em {BoomerAMG}: a parallel algebraic
  multigrid solver and preconditioner}, Applied Numerical Mathematics, 41
  (2000), pp.~155--177.

\bibitem{KimSPE11}
{\sc J.~Kim, H.~A. Tchelepi, and R.~Juanes}, {\em Stability, accuracy, and
  efficiency of sequential methods for coupled flow and geomechanics}, {SPE}
  Journal, 16 (2011), pp.~249--262, \url{https://doi.org/10.2118/119084-pa},
  \url{https://doi.org/10.2118/119084-pa}.

\bibitem{KimTchJua13}
{\sc J.~Kim, H.~A. Tchelepi, and R.~Juanes}, {\em {Rigorous Coupling of
  Geomechanics and Multiphase Flow with Strong Capillarity}}, SPE J., 18
  (2013), pp.~1123--1139, \url{https://doi.org/10.2118/141268-PA}.

\bibitem{Lacroix03}
{\sc S.~Lacroix, Y.~Vassilevski, J.~Wheeler, and M.~F. Wheeler}, {\em Iterative
  solution methods for modeling multiphase flow in porous media fully
  implicitly}, {SIAM} Journal on Scientific Computing, 25 (2003), pp.~905--926,
  \url{https://doi.org/10.1137/s106482750240443x},
  \url{http://dx.doi.org/10.1137/S106482750240443X}.

\bibitem{Lee17}
{\sc J.~J. Lee, K.-A. Mardal, and R.~Winther}, {\em Parameter-robust
  discretization and preconditioning of biot's consolidation model}, {SIAM}
  Journal on Scientific Computing, 39 (2017), pp.~A1--A24,
  \url{https://doi.org/10.1137/15m1029473},
  \url{https://doi.org/10.1137/15m1029473}.

\bibitem{LewSch98}
{\sc R.~W. Lewis and B.~A. Schrefler}, {\em {The Finite Element Method in the
  Static and Dynamic Deformation and Consolidation of Porous Media}}, Wiley,
  Chichester, UK, 2nd~ed., 1998.

\bibitem{Liu15}
{\sc H.~Liu, K.~Wang, and Z.~Chen}, {\em A family of constrained pressure
  residual preconditioners for parallel reservoir simulations}, Numerical
  Linear Algebra with Applications, 23 (2015), pp.~120--146,
  \url{https://doi.org/10.1002/nla.2017},
  \url{https://doi.org/10.1002/nla.2017}.

\bibitem{Mikeli12}
{\sc A.~Mikeli{\'{c}} and M.~F. Wheeler}, {\em Convergence of iterative
  coupling for coupled flow and geomechanics}, Computational Geosciences, 17
  (2012), pp.~455--461, \url{https://doi.org/10.1007/s10596-012-9318-y},
  \url{https://doi.org/10.1007/s10596-012-9318-y}.

\bibitem{Pre14}
{\sc J.~H. Prevost}, {\em {Two-way coupling in reservoir-geomechanical models:
  vertex-centered Galerkin geomechanical model cell-centered and
  vertex-centered finite volume reservoir models}}, Int. J. Numer. Meth. Eng.,
  98 (2014), pp.~612--624, \url{https://doi.org/10.1002/nme.4657}.

\bibitem{Ries79}
{\sc M.~Ries and U.~Trottenberg}, {\em {MGR-Ein} blitzschneller elliptischer
  {L{\"o}ser}}, Preprint 277SBF 72, Universit{\"a}t Bonn,  (1979).

\bibitem{Ries83}
{\sc M.~Ries, U.~Trottenberg, and G.~Winter}, {\em A note on {MGR} methods},
  Linear Algebra and its Applications, 49 (1983), pp.~1 -- 26,
  \url{https://doi.org/10.1016/0024-3795(83)90091-5},
  \url{http://www.sciencedirect.com/science/article/pii/0024379583900915}.

\bibitem{Saad86}
{\sc Y.~Saad and M.~H. Schultz}, {\em {GMRES}: A generalized minimal residual
  algorithm for solving nonsymmetric linear systems}, {SIAM} Journal on
  Scientific and Statistical Computing, 7 (1986), pp.~856--869,
  \url{https://doi.org/10.1137/0907058},
  \url{http://dx.doi.org/10.1137/0907058}.

\bibitem{Scheichl03}
{\sc R.~Scheichl, R.~Masson, and J.~Wendebourg}, {\em Decoupling and block
  preconditioning for sedimentary basin simulations}, Computational
  Geosciences, 7 (2003), pp.~295--318,
  \url{https://doi.org/10.1023/b:comg.0000005244.61636.4e},
  \url{https://doi.org/10.1023/b:comg.0000005244.61636.4e}.

\bibitem{Settari98}
{\sc A.~Settari and F.~Mourits}, {\em A coupled reservoir and geomechanical
  simulation system}, {SPE} Journal, 3 (1998), pp.~219--226,
  \url{https://doi.org/10.2118/50939-pa},
  \url{https://doi.org/10.2118/50939-pa}.

\bibitem{Set_etal17}
{\sc R.~R. Settgast, P.~Fu, S.~D.~C. Walsh, J.~A. White, C.~Annavarapu, and
  F.~J. Ryerson}, {\em {A fully coupled method for massively parallel
  simulation of hydraulically driven fractures in 3-dimensions}}, International
  Journal for Numerical and Analytical Methods in Geomechanics, 41 (2017),
  pp.~627--653, \url{https://doi.org/10.1002/nag.2557}.

\bibitem{Stueben01}
{\sc K.~St{\"u}ben}, {\em A review of algebraic multigrid}, Journal of
  Computational and Applied Mathematics, 128 (2001), pp.~281 -- 309,
  \url{https://doi.org/10.1016/S0377-0427(00)00516-1},
  \url{http://www.sciencedirect.com/science/article/pii/S0377042700005161}.

\bibitem{Stueben07}
{\sc K.~St{\"u}ben, T.~Clees, H.~Klie, B.~Lu, and M.~F. Wheeler}, {\em
  Algebraic multigrid methods ({AMG}) for the efficient solution of fully
  implicit formulations in reservoir simulation}, in {SPE} Reservoir Simulation
  Symposium, Society of Petroleum Engineers ({SPE}), 2007,
  \url{https://doi.org/10.2118/105832-ms},
  \url{http://dx.doi.org/10.2118/105832-MS}.

\bibitem{Trottenberg00}
{\sc U.~Trottenberg, C.~W. Oosterlee, and A.~Schuller}, {\em Multigrid},
  Academic Press, Inc., Orlando, FL, USA, 2001.

\bibitem{Wallis83}
{\sc J.~R. Wallis}, {\em Incomplete gaussian elimination as a preconditioning
  for generalized conjugate gradient acceleration}, in {SPE} Reservoir
  Simulation Symposium, Society of Petroleum Engineers, 1983,
  \url{https://doi.org/10.2118/12265-ms},
  \url{https://doi.org/10.2118/12265-ms}.

\bibitem{Wallis85}
{\sc J.~R. Wallis, R.~P. Kendall, and L.~E. Little}, {\em Constrained residual
  acceleration of conjugate residual methods}, in {SPE} Reservoir Simulation
  Symposium, Society of Petroleum Engineers ({SPE}), 1985.

\bibitem{Wang17}
{\sc L.~Wang, D.~Osei-Kuffuor, R.~D. Falgout, I.~D. Mishev, and J.~Li}, {\em
  Multigrid reduction for coupled flow problems with application to reservoir
  simulation}, in {SPE} Reservoir Simulation Conference, Society of Petroleum
  Engineers ({SPE}), {SPE-182723-MS}, 2017.

\bibitem{White11}
{\sc J.~A. White and R.~I. Borja}, {\em Block-preconditioned {Newton--Krylov}
  solvers for fully coupled flow and geomechanics}, Computational Geosciences,
  15 (2011), pp.~647--659, \url{https://doi.org/10.1007/s10596-011-9233-7},
  \url{http://dx.doi.org/10.1007/s10596-011-9233-7}.

\bibitem{White18}
{\sc J.~A. {White}, N.~{Castelletto}, S.~{Klevtsov}, Q.~M. {Bui},
  D.~{Osei-Kuffuor}, and H.~A. {Tchelepi}}, {\em {A Two-Stage Preconditioner
  for Multiphase Poromechanics in Reservoir Simulation}}, arXiv e-prints,
  (2018), arXiv:1812.05540, p.~arXiv:1812.05540,
  \url{https://arxiv.org/abs/1812.05540}.

\bibitem{White16}
{\sc J.~A. White, N.~Castelletto, and H.~A. Tchelepi}, {\em Block-partitioned
  solvers for coupled poromechanics: A unified framework}, Computer Methods in
  Applied Mechanics and Engineering, 303 (2016), pp.~55--74,
  \url{https://doi.org/10.1016/j.cma.2016.01.008},
  \url{https://doi.org/10.1016/j.cma.2016.01.008}.

\bibitem{Zhou12}
{\sc Y.~Zhou, Y.~J., and H.~A. Tchelepi}, {\em A scalable multistage linear
  solver for reservoir models with multisegment wells}, Computational
  Geosciences, 17 (2012), pp.~197--216,
  \url{https://doi.org/10.1007/s10596-012-9324-0},
  \url{https://doi.org/10.1007/s10596-012-9324-0}.

\end{thebibliography}

\end{document}